\documentclass[lineno]{biometrika}

\usepackage{times,url}
\usepackage{bm}
\usepackage{natbib,amsmath,epsfig,epsf,psfrag,ntheorem,amssymb,psfig}
\RequirePackage{fix-cm}

\usepackage[plain,noend]{algorithm2e}

\makeatletter
\renewcommand{\algocf@captiontext}[2]{#1\algocf@typo. \AlCapFnt{}#2} 
\def\@algocf@capt@plain{top}
\renewcommand{\algocf@makecaption}[2]{%
  \addtolength{\hsize}{\algomargin}%
  \sbox\@tempboxa{\algocf@captiontext{#1}{#2}}%
  \ifdim\wd\@tempboxa >\hsize
    \hskip .5\algomargin%
    \parbox[t]{\hsize}{\algocf@captiontext{#1}{#2}}
  \else%
    \global\@minipagefalse%
    \hbox to\hsize{\box\@tempboxa}
  \fi%
  \addtolength{\hsize}{-\algomargin}%
}
\makeatother

\def\T{{ \mathrm{\scriptscriptstyle T} }}

\begin{document}

\newcommand{\wt}{\widetilde}
\newcommand{\wh}{\widehat}
\newcommand{\hatb}{\hat b}
\newcommand{\ol}{\overline}



\markboth{X. Wang, P. Du and J. Shen}{Spatially adaptive smoothing splines}

\title{Smoothing splines with varying smoothing parameter}

\author{XIAO WANG}
\affil{Department of Statistics, Purdue University, 250 N.
University Street, West Lafayette, Indiana, 47907, USA
\email{wangxiao@purdue.edu} }

\author{PANG DU}
\affil{Department of Statistics, Virginia Tech, 406-A Hutcheson Hall
Blacksburg, Virginia 24061, USA \email{pangdu@vt.edu}}

\author{\and JINGLAI SHEN}
\affil{Department of Mathematics and Statistics, University of
Maryland at Baltimore County, 1000 Hilltop Circle, Baltimore,
Maryland, 21250, USA \email{shenj@umbc.edu}}

\maketitle

\begin{abstract}
This paper considers the development of spatially adaptive smoothing
splines for the estimation of a regression function with non-homogeneous smoothness across the domain.
Two challenging
issues that arise in this context are the evaluation of the equivalent kernel and the determination
of a local penalty.
The roughness penalty is a function
of the design points in order to accommodate local behavior of the
regression function.
It is shown that the spatially adaptive smoothing
spline estimator is approximately a kernel estimator. The
resulting equivalent kernel is spatially dependent. The equivalent kernels
for traditional smoothing splines are a special case of this general
solution. With the aid of the Green's function for a two-point boundary value problem, the explicit forms of the
asymptotic mean and variance are obtained for any interior point. Thus, the optimal
roughness penalty function is obtained by approximately minimizing the
asymptotic integrated mean square error.
Simulation results and an application illustrate the performance of the proposed estimator.
\end{abstract}


\section{Introduction}

Smoothing splines play a central role in
nonparametric curve-fitting. Recent synopses include
\cite{wahba_90}, \cite{eubank_99}, \cite{gu_02}, and
\cite{eggermont_09}. Specifically, consider the problem of
estimating the mean function from a regression model
$$y_i = f_0(t_i) + \sigma(t_i)\epsilon_i ~ (i=1, \ldots, n),$$
where the $t_i$ are the design points on $[0, 1]$,
the $\epsilon_i$ are independent and  identically-distributed random
variables with zero mean and unit variance, $\sigma^2(\cdot)$ is
the variance function, and $f_0$ is the underlying true regression
function.  The traditional smoothing spline is formulated as the
solution $f$ to the minimization of
\begin{equation}\label{equ:ss}
{1\over n}\sum_{i=1}^n\sigma^{-2}(t_i)\{y_i - f(t_i)\}^2 + \lambda \int_0^1 \{f^{(m)}(t)\}^2dx,
\end{equation}
where $\lambda>0$ is the penalty parameter controlling the
trade-off between the goodness-of-fit and smoothness of the fitted
function. Smoothing splines have a solid theoretical foundation
and are among the most widely used methods for nonparametric
regression \citep{speckman_81, cox_83}.

The traditional smoothing spline model has a major deficiency:
it uses a global smoothing parameter $\lambda$,  so the degree of smoothness of $f_0$ remains about the
same across the design points. This makes it
difficult to efficiently estimate functions with non-homogeneous
smoothness. \cite{wahba_95} suggested
using a more general penalty term, which replaces the constant
$\lambda$ by a roughness penalty function $\lambda(\cdot)$. Since
$\lambda(\cdot)$ is then a function of $t$, the model
becomes adaptive in the sense that it accommodates the local
behavior of $f_0$ and imposes a heavier penalty in the regions of
lower curvature of $f_0$. \citet{pintore_06} used a piecewise
constant approximation for $\lambda(\cdot)$ but this requires specification of the number of knots, the knot
locations, and the values of $\lambda(\cdot)$ between these
locations. \cite{Storlie_09JCGS} discussed some
computational issues on spatially adaptive smoothing splines.
\cite{guo_10} refined the piecewise constant idea and designed a data-driven algorithm to determine the optimal jump locations and sizes for $\lambda(\cdot)$. Besides adaptive smoothing splines,
other adaptive methods have been developed, including
variable-bandwidth kernel smoothing (M$\mathrm{\ddot{u}}$ller \&\
Stadtm$\mathrm{\ddot{u}}$ller, 1987), adaptive wavelet shrinkage
\citep{donoho_94,donoho_95,donoho_98}, local polynomials with
variable bandwidth \citep{fan_96}, local penalized splines
\citep{ruppert_00}, regression splines \citep{friedman_89,
stone_97, luo_97, hansen_02}, and free-knot splines
\citep{mao_03}. Further, Bayesian adaptive regression has also
been reported by \cite{smith_96}, \cite{DiMatteo_01}, and
\cite{wood_02}. Nevertheless, adaptive smoothing splines have the
advantages of computational efficiency and
easy extension to multidimensional covariates using the smoothing
spline analysis of variance technique \citep{wahba_90, gu_02}. Further, the results in the present paper can be extended to the
more general L-spline smoothing \citep{kim_71, kohn_83, wahba_85}.
Also, the usual Reinsch scheme can be easily modified to the present case.

Let
$
 W^m_2 = \{f: f^{(m-1)} \mbox{
absolutely continuous and } f^{(m)}\in L_2[0,1]\}$,  where $L_2[0,
1]$ is the space of Lebesgue square integrable functions,
endowed with its usual norm $\|\cdot\|_2$ and inner product $(
\cdot,~ \cdot )_2$. The method of adaptive smoothing splines
finds $f \in W^m_2$ to minimize the functional
\begin{equation}\label{equ:obj}
\psi(f) = {1\over n}\sum_{i=1}^n\sigma^{-2}(t_i)\{y_i - f(t_i)\}^2 + \lambda \int_0^1 \rho(t) \{f^{(m)}(t)\}^2dt,
\end{equation}
where $\lambda>0$ is the penalty parameter, and $\rho:[0, 1]
\rightarrow \mathbb (0 , \infty)$ denotes the adaptive penalty
function; more properties of $\rho$ will be stated later. Here, we
incorporate a function $\rho(t)$ into the roughness penalty,
which generalizes the traditional smoothing splines, where
$\rho(t)\equiv 1$.  
A two-point boundary value problem technique
has been developed
to find the asymptotic mean squared error of the adaptive
smoothing spline estimator with the aid of the Green's function.
Thus the optimal roughness penalty function is obtained
explicitly by approximately minimizing the asymptotic integrated mean squared error. Asymptotic analysis of traditional smoothing
splines using Green's functions was performed by \cite{rice_83},
\cite{silverman_84}, \cite{messer_91}, \cite{nychka_95}, and \cite{eggermont_09}; an
extension to certain adaptive splines was made in
\cite{abramovich_99}. In contrast to these results, the
current paper develops a general framework for asymptotic analysis
of adaptive smoothing splines, yielding a
systematic, yet relatively simpler, approach to obtaining closed-form
expressions of equivalent kernels for interior points and to
asymptotic analysis. Our estimate possesses the interpretation of spatial adaptivity \citep{donoho_98}, and the equivalent kernel may vary in shape and bandwidth from point to point,
depending on the data.



\section{Characterizations of the estimator}

In this section, we derive the optimality conditions for
 the solution that minimizes the functional
(\ref{equ:obj}). Let
$\omega_n(t) = n^{-1}\sum_{i=1}^n {\cal I}(t_i\le t)$
where ${\cal I}$ is the indicator function, and let $\omega$ be a distribution
function with a continuous and strictly positive density function
$q$ on $[0, 1]$. For a function $g$, define $\|g\| = \sup_{t\in[0,
1]}|g(t)|$ and subsequent norms likewise.  Let $D_n = \|\omega_n - \omega\|$.
If the design points $t_i$ are equally spaced, $D_n = O(n^{-1})$ with $q(t)=1$ for $t\in [0,1]$. If $t_i$ are
 independent and identically distributed regressors from a
distribution with bounded positive density $q$, then $D_n =
O\{n^{-1/2}(\log\log n)^{1/2}\}$ by the law of the iterated
logarithm for empirical distribution functions.

Let $h$ be a piecewise constant function such that
$h(t_i)=y_i$ ($i=1, \ldots, n$). For any $t\in [0, 1]$ and $f\in L_1[0, 1]$,
define
$$l_1(f, t) = \int_0^t \sigma^{-2}(s)f(s)d\omega(s),~~~~ ~l_k(f, t) = \int_0^t l_{k-1}(f, s)ds,$$
and
$$\check l_1(f, t) = \int_0^t \sigma^{-2}(s)f(s)d\omega_n(s),~ ~~~ ~ \check l_k(f, t) = \int_0^t \check l_{k-1}(f, s)ds~~~~ (2\le k\le m).$$


\begin{theorem}\label{thm:char} Necessary and sufficient conditions for $\hat f\in W^m_2$ to minimize $\psi$ in (\ref{equ:obj}) are that
\begin{equation}\label{equ:optima_condition1}
(-1)^m~ \lambda ~\rho(t)~ \hat f^{(m)}(t) + \check l_m(\hat f, t) = \check l_m(h, t),  t\in [0,1],\end{equation}
almost everywhere, and
\begin{equation}\label{equ:optima_condition2}\check l_k(\hat f, 1) = \check l_k(h, 1) ~~(k=1, \ldots, m).\end{equation}
\end{theorem}

Both $\check l_1(\hat f, t)$ and $\check l_1(h, t)$  are
piecewise constant in $t$. Therefore $\check l_m(h, t) - \check
l_m(\hat f, t)$ is a piecewise $(m-1)$th order polynomial. Thus,
 Theorem \ref{thm:char} shows that $\rho(t)~ \hat
f^{(m)}(t)$ is a piecewise $(m-1)$th order polynomial. The exact
form of $\hat f$ will depend on additional assumptions about
$\rho(t)$. For example, \cite{pintore_06} assumed $\rho(t)$ to be
piecewise-constant with possible jumps at a subset of the design
points. Then, the optimal solution is a polynomial spline of order
$2m$. It is well-known that the traditional smoothing spline is a
natural spline of order $2m$, which corresponds to the case here when
$\rho(t)\equiv 1$.

\section{Asymptotic properties of the estimator}\label{sec:asy}

We establish an equivalent kernel and asymptotic
distribution of the spatially adaptive smoothing splines at interior points using a two-point boundary value problem technique. The key idea is to represent the solution to
(\ref{equ:optima_condition1}) by a Green's function. It will be
shown that the adaptive smoothing spline estimator can be approximated by a
kernel estimator, using this Green's function.

Denote $R_k(t) = l_k(\hat f, t) - \check l_k(\hat f, t)$ ($k=1,
\ldots, m$).  Specifically, when $k=m$, it follows from Theorem
\ref{thm:char} that
$$R_m(t) = (-1)^m \lambda~ \rho(t) \hat f^{(m)}(t) +  l_m(\hat f, t) - l_m(h, t).$$
Write $r(t) = \sigma^2(t)/q(t)$. Thus, $l_m(\hat f, t)$ solves the two-point boundary value problem
\begin{equation}\label{equ:ode_rep}
(-1)^m \lambda~ \rho(t) {d^m\over dt^m}\Big\{  r(t) {d^m\over dt^m}l_m(\hat f, t)\Big\} + l_m(\hat f, t) = \check l_m(h, t)+ R_m(t),\end{equation}
 subject to the $2m$ boundary conditions from (\ref{equ:optima_condition2}):
\begin{equation}\label{equ:ode_rep_boundary}
 l_k(\hat f, 0)=0, l_k(\hat f, 1) = l_k(h, 1)+R_k(1) (k=1, \ldots, m).
  \end{equation}

The solution to (\ref{equ:ode_rep}) can be obtained explicitly with the aid of the Green's function.
For readers unfamiliar with
Green's functions, operationally speaking, if $P(t, s)$ is the
Green's function for
\begin{equation}\label{equ:ode00}
(-1)^m \lambda \rho(t) \{r(t)
u^{(m)}(t)\}^{(m)} + u(t) = 0,\end{equation}
then $\int_0^1 P(t,
s)\{\check l_m(h, s) + R_m(s)\}ds$ will solve (\ref{equ:ode_rep}).
This, together with the
boundary conditions (\ref{equ:ode_rep_boundary}),
yields the solution to the two-point boundary value problem in (\ref{equ:ode_rep}) and (\ref{equ:ode_rep_boundary}). The derivations of the Green's function and discussions of the boundary conditions are given in the online Supplementary Material. Specifically, let $\{C_{k}(t), k=1, \ldots, 2m\}$ be $2m$ linearly independent solutions for the homogeneous differential equation
$$(-1)^m \lambda~ \rho(t) {d^m\over dt^m}\Big\{ r(t){d^m\over dt^m}l_m(\hat f, t)\Big\} + l_m(\hat f, t)  =0.$$
Then, $l_m(\hat f, t)$ in (\ref{equ:ode_rep}) can be represented as
\begin{equation}\label{equ:repF}
l_m(\hat f, t) =\int_0^1 P(t, s)\check l_m(h, s)ds + \int_0^1 P(t, s)R_m(s)ds + \sum_{k=1}^{2m} a_k C_k(t),
\end{equation}
where the last term is due to the boundary conditions and the
coefficients $a_k  (k=1, \ldots, 2m)$ are shown to be unique
and stochastically bounded for all sufficiently small $\lambda$ in
the  Supplementary Material. Equation (\ref{equ:repF}) can be
decomposed into three parts:  the asymptotic mean $\int_0^1
P(t, s) l_m(f_0, s)ds$; the random component $\int_0^1 P(t,
s) \check l_m(h-f_0, s)ds$; and the
remainder term $\Gamma(t)= \sum_{k=1}^{2m} a_k C_k(t)+\int_0^1 P(t, s) \tilde R_m(s) ds$,
where $\tilde R_m(t) = l_m(\hat f-f_0, t) - \check l_m(\hat f-
f_0, t)$.
It will be shown that $\|\tilde R_m\|$ has a smaller order and the
remainder term is negligible in the asymptotic analysis.
Taking the $m$-th derivative point-wise on both sides of
(\ref{equ:repF}) gives the crucial representation of the adaptive
smoothing spline estimator. This gives
\begin{equation}\label{equ:rep0}
r^{-1}(t)\hat f(t) =  {d^m\over dt^m}\int_0^1 P(t, s) l_m(f_0, s)ds +   {d^m\over dt^m}\int_0^1  P(t, s) \check l_m(h - f_0, s)ds + \Gamma^{(m)}(t).
\end{equation}


We now introduce the main assumptions of this paper:
\begin{assumption}\label{A1}
The functions $\rho(\cdot)$, $q(\cdot)$, and $\sigma(\cdot)$ are $(m+1)$-times continuously differentiable and strictly positive.
\end{assumption}
\begin{assumption}\label{A2}
The function $f_0$ is $2m$-times continuously differentiable.
\end{assumption}
\begin{assumption}\label{A3}
The smoothing parameter
$\lambda\rightarrow 0$ as $n\rightarrow \infty$.
Denote
$$\Delta_n = D_n n^{-1/2} \lambda^{-(1+m)/(2m)}{ \max\Big[ \{\log (1/ \lambda)\}^{1/2}, (\log\log n)^{1/2}\Big] } .$$
Assume $\Delta_n\rightarrow 0$ as $n\rightarrow\infty$.
\end{assumption}
\begin{assumption}\label{A4}
The random errors $\epsilon_i$ have a
finite fourth moment.
\end{assumption}


Assumption \ref{A3} ensures that the smoothing
parameter $\lambda$ tends to zero not too quickly. In particular,
it encompasses the cases of equally spaced design variables and of independent and identically-distributed regressors from a
distribution with bounded positive density. In the former case,
$D_n = O(n^{-1})$ and in the second case, $D_n =
O(n^{-1/2}(\log\log n)^{1/2})$. The optimal choice
of $\lambda$ discussed subsequently is of order $n^{-2m/(4m+1)}$
and it is easy to check that it satisfies Assumption \ref{A3}.

\begin{theorem}\label{thm:asy-rep} Assume that Assumptions \ref{A1}--\ref{A4} hold.  Let $\beta = \lambda^{-{1/(2m)}}$. For any given $t\in
(0, 1)$, the adaptive smoothing spline estimator $\hat f$ can be
written as
\begin{eqnarray}\label{equ:repf}
\hat f(t) &=& f_0(t) + \lambda~ (-1)^{m-1} r(t) \big\{\rho(t) f_0^{(m)}(t)\big\}^{(m)} + o(\lambda) + {1\over n} \sum_{i=1}^n {\sigma(t_i)\over q(t_i)} J(t, t_i)\epsilon_i\\
&& + O(\beta^m)\Delta_n + O(\beta^m)e^{-\beta O(1)} \nonumber
\end{eqnarray}
uniformly in $\lambda$, where $J(t, s)$ is given in (\ref{equ:J}).
\end{theorem}

\begin{remark}\rm  \cite{eggermont_06} were the first to show in full generality that the
standard spline smoothing corresponds approximately to smoothing
by a kernel method. A simple explicit formula of the equivalent kernel for all $m$, denoted by $K(t, s)$, is given by \cite{berlinet_04}.
For interior points, the kernel $K$ is of the form $K(t, s)
= \beta L(\beta |t-s|)$ for some function $L$, and $L(|\cdot|)$ is
a $2m$-th order kernel on $(-\infty, \infty)$. In particular, the
shape of $K(t, \cdot)$ is defined by $L(\cdot)$ and is the same
for different $t$. For example, the closed form expressions for the
first two equivalent kernels are:
\begin{eqnarray*}
   m=1: \ \ L(|t|)  & = & \frac{1}{2} e^{-|t|}, \\
   m=2: \ \ L(|t|)  & = & \frac{1}{2^{3/2}} \, e^{- |t|/2^{1/2}}
    \Big\{ \, \cos\Big(\frac{|t|}{2^{1/2}}\Big)+ \sin\Big( \frac{|t|}{2^{1/2}}\Big) \, \Big\},\\
   m=3: \ \ L(|t|)  & = & \frac{1}{6} e^{-|t|} +  \, e^{-\frac{1}{2} |t|}
    \Big\{ \, \frac{1}{6} \cos\Big(\frac{3^{1/2}|t|}{2}\Big) + \frac{3^{1/2}}{6} \sin\Big(\frac{3^{1/2}|t|}{2}\Big)
     \, \Big\}, \\
   m=4: \ \ L(|t|)  & = &  e^{-0{\cdot}9239 |t|}
    \Big\{ \, 0{\cdot}2310 \cos(0{\cdot}3827|t|) + 0{\cdot}0957 \sin(0{\cdot}3827|t|)
     \, \Big\}  \\
   & &  \, + \, e^{-0{\cdot}3827 |t|} \Big\{ \, 0{\cdot}0957  \cos(0{\cdot}9239|t|) + 0{\cdot}2310 \sin(0{\cdot}9239|t|)
     \, \Big\}.
\end{eqnarray*}
Theorem \ref{thm:asy-rep} indicates that the spatially adaptive smoothing spline estimator is also
approximately a kernel regression estimator. The equivalent kernel
$J(t, s)$ is the corresponding Green's function. As shown in the Supplementary Material,
\begin{equation}\label{equ:J}
J(t, s) = \beta \varrho(s) Q_\beta'(s) L\{\beta|Q_\beta(t)-Q_\beta(s)|\},
\end{equation}
where
$$Q_\beta(t) = \int_0^t  \Big\{r(s) \rho(s)\Big\}^{-{1/(2m)}}\Big\{1+O(\beta^{-1})\Big\}ds$$
is an increasing function of $t$, and
$\|\varrho\|=1+O(\beta^{-1})$. This shows that the shape of $J(t, \cdot)$ varies with $t$. Our estimator possesses the interpretation of spatial adaptivity \citep{donoho_98}; it is
asymptotically equivalent to a kernel estimator with a kernel that varies in shape and bandwidth from point to point.
\end{remark}


\begin{remark}\rm
The number $\beta^{-1}$ in (\ref{equ:J}) plays a role similar to the
bandwidth $h$ in kernel smoothing. Theorem
\ref{thm:asy-rep} shows that the asymptotic mean has bias $
(-1)^{m-1} \lambda r(t) \big\{\rho(t)
f_0^{(m)}(t)\big\}^{(m)}$, which can be negligible if $\lambda$ is
reasonably small. On the other hand, $\lambda$ cannot be
arbitrarily small since that will inflate the random component. The
admissible range for $\lambda$ is a compromise between these two.
\end{remark}



\begin{corollary}\label{coro:asy} Given $\rho(\cdot)$ and $r(\cdot)$, and assuming Assumptions \ref{A1}--\ref{A4}, if $\lambda = n^{-2m/(4m+1)}$, then, for any $t\in (0,1)$, $n^{2m/ (4m+1)}\big\{\hat f(t) - f_0(t)\big\}$ converges to
\begin{equation}
 N\Big[(-1)^{m-1} r(t) \big\{\rho(t) f_0^{(m)}(t)\big\}^{(m)}, ~~~~L_0~r(t)^{1-{1/ (2m)}} \rho(t)^{-{1/ (2m)}}   \Big],
\end{equation}
in distribution, where $L_0 = \int_{-\infty}^\infty L^2(|t|)dt$.
\end{corollary}

The proof of Corollary \ref{coro:asy} is given in the Supplementary Material.
The asymptotic mean
squared error of the spatially adaptive smoothing spline estimator is of
order $n^{-4m/(4m+1)}$, which is the optimal rate of
convergence given in \cite{stone_82}.

\section{Optimal selection of $\rho$}

The optimal $\lambda$ and $\rho$ are chosen to minimize the integrated asymptotic mean squared error
\begin{equation}\label{equ:mse}\int_0^1 \Big\{\lambda^2 r^2(t)\big[\{\rho(t) f_0^{(m)}(t)\}^{(m)}\big]^2 + {L_0\over n\lambda^{1/(2m)}}r(t)^{1-{1/ (2m)}} ~\rho(t)^{-{1/ (2m)}}\Big\}dt,\end{equation}
which is in fact a function of $\lambda \rho(t)$.
We choose the optimal $\lambda$ to be $\lambda^o = n^{-2m/(4m+1)}$.
The optimal roughness penalty function $\rho(t)$ minimizes the functional
\begin{equation}\label{equ:mse2}\Pi(\rho) = \int_0^1 \left\{r^2(t)\big[\{\rho(t) f_0^{(m)}(t)\}^{(m)}\big]^2 + L_0 r(t)^{1-{1/ (2m)}} ~\rho(t)^{-{1/ (2m)}}\right\}dt.\end{equation}

  Without any further assumptions, the above minimization problem  does not have an optimal solution, since any arbitrarily large and positive function $\rho$ with $\big\{\rho(t) f_0^{(m)}(t) \big\}^{(m)} =0$ on any sub-interval of $[0, 1]$ will make $\Pi(\cdot)$ arbitrarily small. To deal with this problem, we first impose a technical assumption on $f_0$.
\begin{assumption}\label{A5}
The set $\mathcal N=\big\{ t \in [0, 1]: f^{(m)}_0(t)=0 \big\}$ has zero measure.
\end{assumption}

Let $u(t)= \{\rho(t) f_0^{(m)}(t)\}^{(m)}$, $z(t)= \rho(t) f^{(m)}_0(t)$, and $D^{-m}$ be the $m$-fold integral operator. Then $z^{(m)}(t) = u(t)$ and
\begin{equation} \label{eqn:z}
    z(t) = (D^{-m} u)(t) + {\theta}^\T(t) x_0,
\end{equation}
for $\theta(t) = \big(1, t, t^2/2!, \ldots,  t^{m-1}/ (m-1)!\big)^\T$ and some $x_0 \in \mathbb R^m$.
Moreover, we can define $z(t)/f^{(m)}_0(t)$ to be any positive constant
for all $t\in \mathcal N$ where $f^{(m)}_0(t)=0$. This definition is assumed in the subsequent development. Hence, the functional $\Pi(\rho)$ in (\ref{equ:mse2}) becomes
\[
    J(u, x_0)  = \int_0^1 r^2(t) u^2(t) dt + \int^1_0 L_0 r(t)^{1-{1/(2m)}} ~\left\{ { z(t) \over f^{(m)}_0(t) } \right\}^{-{1/ (2m)}} dt,
\]
where $z(t)$ is defined by $(u, x_0)$. We then introduce another technical assumption on $z(t)$, or essentially on $\rho$.
\begin{assumption}\label{A6}
There exist positive constants $\mu$ and $\varepsilon$ such that $\| x_0\| \le \mu$ and $z(t)/f^{(m)}_0(t)\ge \varepsilon$ for all $t$. And $\big\{ z(t)/f^{(m)}_0(t)  \big\}^{-1/(2m)}$ is Lebesgue integrable on  $[0, 1]$.
\end{assumption}


Consider the following set in $L_2[0, 1] \times \mathbb R^m$,
\begin{multline*}
  \mathcal P  =  \Big\{ (u, x_0) \in L_2[0, 1] \times \mathbb R^m \, : \, \| x_0 \| \le \mu, \ z(t)/f^{(m)}_0(t) \ge \varepsilon \ \mbox{ for all } t \in [0, 1], \  \mbox{ and } \\ \Big.
   \Big.  \big\{ z(t)/f^{(m)}_0(t)  \big\}^{-1/(2m)}  \mbox{ is Lebesgue integrable on }    [0, 1] \,  \Big\},
\end{multline*}
where $z(t)$ is given in (\ref{eqn:z}) dependent on $(u, x_0)$.
Further development in the Supplemental Material establishes the following theorem that the objective functional $J$ attains a unique minimum in $\mathcal{P}$.  In fact, under the additional Assumptions \ref{A5} and \ref{A6}, the theorem first shows the existence of an optimal solution. Moreover, since the objective functional $J$ is strictly
convex and the constraint set $\mathcal{P}$ is convex, the uniqueness of an optimal solution also follows.

%
%

\begin{theorem} \label{thm:existene_unique}
  Under Assumptions \ref{A1}, \ref{A2}, \ref{A5} and \ref{A6}, the optimization problem
  $
            \inf_{(u, x_0)\in \mathcal P} J(u, x_0)
  $
  has a unique solution in $\mathcal P$.
\end{theorem}

\begin{remark}
  Given the optimal solution $(u^*, x^*)$, $z_{(u^*, x^*)}(t)$  is bounded on $[0, 1]$ due to its absolute continuity. The lower bound $\varepsilon$ in Assumption \ref{A6} ensures that the optimal $\rho$ is bounded below from zero. However, there is no guarantee that the optimal $\rho$ is bounded above due to the possibility for small values of $\big|f^{(m)}_0\big|$.  To avoid this problem, one may impose an additional upper bound constraint in Assumption \ref{A6}. The proof of existence and uniqueness remains the same.
\end{remark}

\section{Implementation}

Obtaining an explicit solution of (\ref{equ:mse2}) is difficult. Motivated from \cite{pintore_06}, we consider approximating $\rho$ by a piecewise constant function such that $\rho(t) = \rho_j$ for $t\in (\tau_{j-1},
\tau_j]$, $j=0, \ldots, S+1$. Here $\tau_0=0, \tau_{S+1}=1$, and
$0<\tau_1<\cdots<\tau_S<1$ are interior adaptive smoothing knots
whose selection will be described below.
When the integral in (\ref{equ:mse2}) is taken ignoring the non-differentiability at the jump points $\tau_j$ $(j=1, \ldots, S)$, we obtain
$$\sum_{j=1}^{S+1} \left[\rho_j^2 \int_{\tau_{j-1}}^{\tau_j} r^2(t) \{f_0^{(2m)}(t)\}^2dt
+ \rho_j^{-{1/ (2m)}} L_0 \int_{\tau_{j-1}}^{\tau_j}
r(t)^{1-{1/ (2m)}}dt\right].$$
Therefore, the optimal $\rho_j$ is
\begin{equation}\label{equ:est-rho}\rho_j = \left[L_0 \int_{\tau_{j-1}}^{\tau_j} r(t)^{1-{1/ (2m)}}dt
\over 4m\int_{\tau_{j-1}}^{\tau_j} r^2(t) \{f_0^{(2m)}(t)\}^2dt\right]^{2m/ (4m+1)}, ~~~~~~~j=1, \ldots, S+1.\end{equation}
Unfortunately, the optimal values for the $\rho_j$ depend on
$r(t)$  and the $2m$-th derivative of the underlying regression function
$f_0(t)$. We replace them by estimates in practice.

\begin{remark}\rm
Rigorously speaking, such a step-function approximation to $\rho$ is not a valid solution to \eqref{equ:mse2} due to non-differentiability.
However, simulations seem to suggest that such a simple approximation can yield good results.
Furthermore, one can modify such $\rho$, for example, to make it satisfy Assumption \ref{A2}. In a sufficiently small neighborhood of each jump point, one can replace the steps by a smooth curve connecting the two steps such that the resulting function satisfies Assumption \ref{A2}. Hence the piecewise constant $\rho$ can be viewed as a simple approximation to this smooth version of $\rho$.
\end{remark}

We now describe the detailed steps for approximate implementation.
The first step is to select the interior smoothing knots
$\tau_{j}$ $(j=1,\ldots,S)$. An abrupt change in the
smoothness of the function is often associated with a similar change
in the conditional probability density of $y$ given $t$. For example, a steeper
part of the function often comes with sparser data, or smaller
conditional probability densities of $y$ given $t$. Hence, we first use the
\texttt{sscden} function in the R package \texttt{gss} to estimate
the conditional probability densities of $y$ given $t$ on a dense grid, say
$s_k=k/100$ $(k=1,\ldots,100)$. Then with a given $S$, we select the
top $S$ $s_k$ where the conditional probability density changes the most
from $s_k$ to $s_{k+1}$. A more accurate but considerably more
time-consuming way of selecting the smoothing knots is a binary
tree search algorithm proposed in \cite{guo_10}.

Estimation of $\sigma^2(t)$ was first studied by \cite{muller_87}. In this paper, we use the local polynomial approach in \cite{fan_98};
see \cite{hall_89}, \cite{ruppert_97}, and \cite{cai_08} for other methods.
This provides the weights for obtaining a weighted smoothing spline estimate
of $f(t)$, whose derivative yields an estimate of $f^{(2m)}(t)$.
The function $q(t)$ can be replaced by an estimate of the density function of $t_i$ $(i=1,\ldots,n)$.
All these computations can be conveniently carried out using the R packages
\texttt{locpol} and \texttt{gss}.

Ideally, the optimal $\rho_j$ computed as above work well.
However, similar to the finding in \cite{Storlie_09JCGS}, we have found that
a powered-up version $\rho_j^\gamma$ for some $\gamma>1$ can often help in practice.
Intuitively, this power-up makes up a bit
for the under-estimated differences in $f^{(2m)}(t)$ across the predictor domain.

For the tuning parameters $S$ and $\gamma$, we consider
$S\in \{0,2,4,8\}$ and $\gamma\in\{1,2,4\}$. Theoretically a larger $S$ might be preferred due to the better approximation of such step functions to the real function. However, as shown in \cite{pintore_06} and \cite{guo_10}, an $S$ greater than 8 tends to overfit the data. The options for $\gamma$ were suggested in \cite{Storlie_09JCGS}.
In traditional smoothing splines, smoothing parameters are selected by the
generalized cross-validation \citep{craven_79} or the generalized maximum likelihood estimate \citep{wahba_85}. As pointed out in \cite{pintore_06},
a proper criterion for selecting the piecewise constant $\rho(\cdot)$ should
penalize on the number of segments of $\rho$. The generalized Akaike information criterion proposed
in \cite{guo_10} serves this purpose, which is a penalized version of the generalized maximum likelihood estimate
where $S$ is penalized similar to the degrees of
freedom in the conventional Akaike information criterion. In this paper, we will use the generalized Akaike information criterion to select $S$ and $\gamma$.

Once the piecewise constant penalty function $\rho$ is determined, we
compute the corresponding adaptive smoothing spline estimate as follows.
By the representer theorem \citep{wahba_90}, the minimizer of
\eqref{equ:obj} lies in a finite-dimensional space of functions
\begin{equation}\label{equ:fexp}
f(t)=\sum_{i=1}^n c_i K_\rho(t_i,t)+\sum_{j=0}^{m-1} d_j\phi_j(t),
\end{equation}
where $c_i$ and $d_j$ are unknown coefficients,
$\phi_j(t)=t^j/j!$ for $j=0,\ldots,m-1$, and $K_\rho$ is the reproducing kernel
function whose closed form expressions at $(t_i,\cdot)$ with a piecewise-constant $\rho$
are given in Section~2.2 of \cite{pintore_06}.
Plugging \eqref{equ:fexp} into \eqref{equ:obj}, we solve for
$c=(c_1,\ldots,c_n)^\T$ and $d=(d_0,\ldots,d_{m-1})^\T$ by
the Newton--Raphson procedure with a fixed $\lambda$. Here $\lambda$ can
be selected by the generalized cross-validation or the generalized maximum likelihood estimate with the adaptive reproducing kernel function.


\section{Simulations}

This section  compares the
estimation performance of different smoothing spline methods.
For traditional smoothing splines, we used the cubic smoothing splines from the function \texttt{ssanova} in the R package \texttt{gss} and the smoothing parameter was selected by the generalized cross validation score.
For the spatially adaptive smoothing splines in \cite{pintore_06}, we used an equally-spaced five-step penalty function following their implementation and the optimal penalty function was selected to minimize the generalized cross validation function (19) in their paper.
For the Loco-Spline in \cite{Storlie_09JCGS}, we downloaded the authors' original program from the site of the Journal of Computational and Graphical Statistics:
\texttt{http://amstat.tandfonline.com/doi/suppl/10.1198/jcgs.2010.09020/}
\texttt{suppl\_file/r-code.zip}.
For the proposed adaptive smoothing splines, we used $m=1$ and cubic smooth splines to compute the optimal $\rho_j$'s.

\begin{figure}
\centerline{\psfig{file=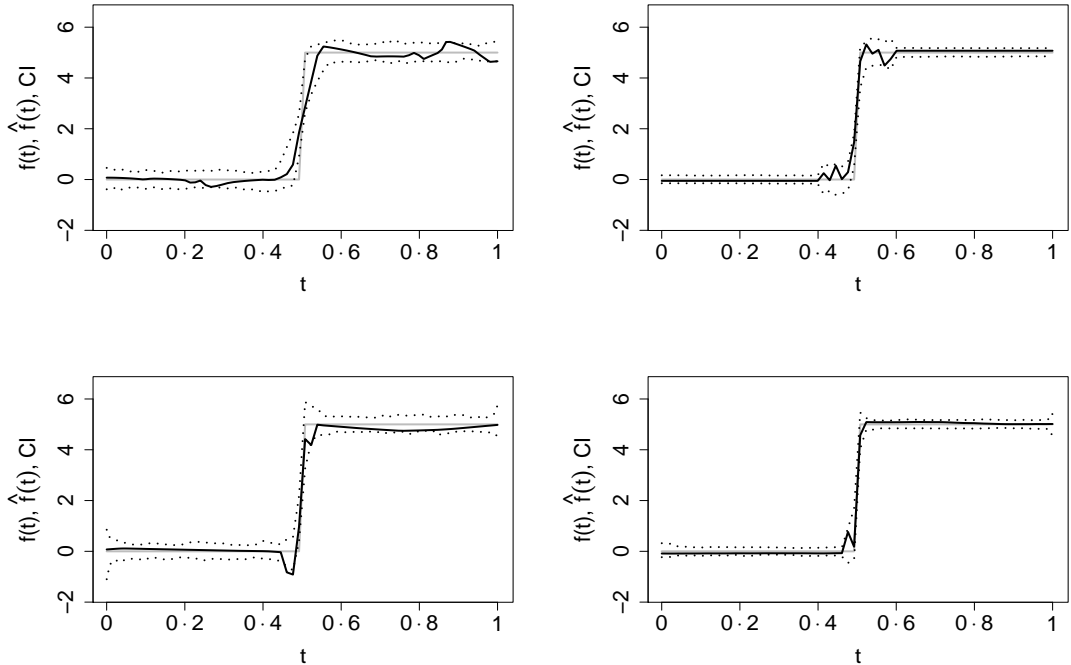,height=\linewidth,width=0.67\linewidth,clip=,angle=270}}
\caption{Estimates of the Heaviside function for the data replicates with median integrated square errors.
The plotted curves are the true function (solid line), the spline estimate (solid line), and the
point-wise empirical 0.025 and 0.975 quantiles (dotted lines).
Top left: traditional smoothing spline estimate. Top right: estimate from the method in \cite{pintore_06}.
Bottom left: Loco-Spline estimate.
Bottom right: proposed adaptive smoothing spline estimate.
}
\label{fig:heav}
\end{figure}

\begin{figure}
\centerline{\psfig{file=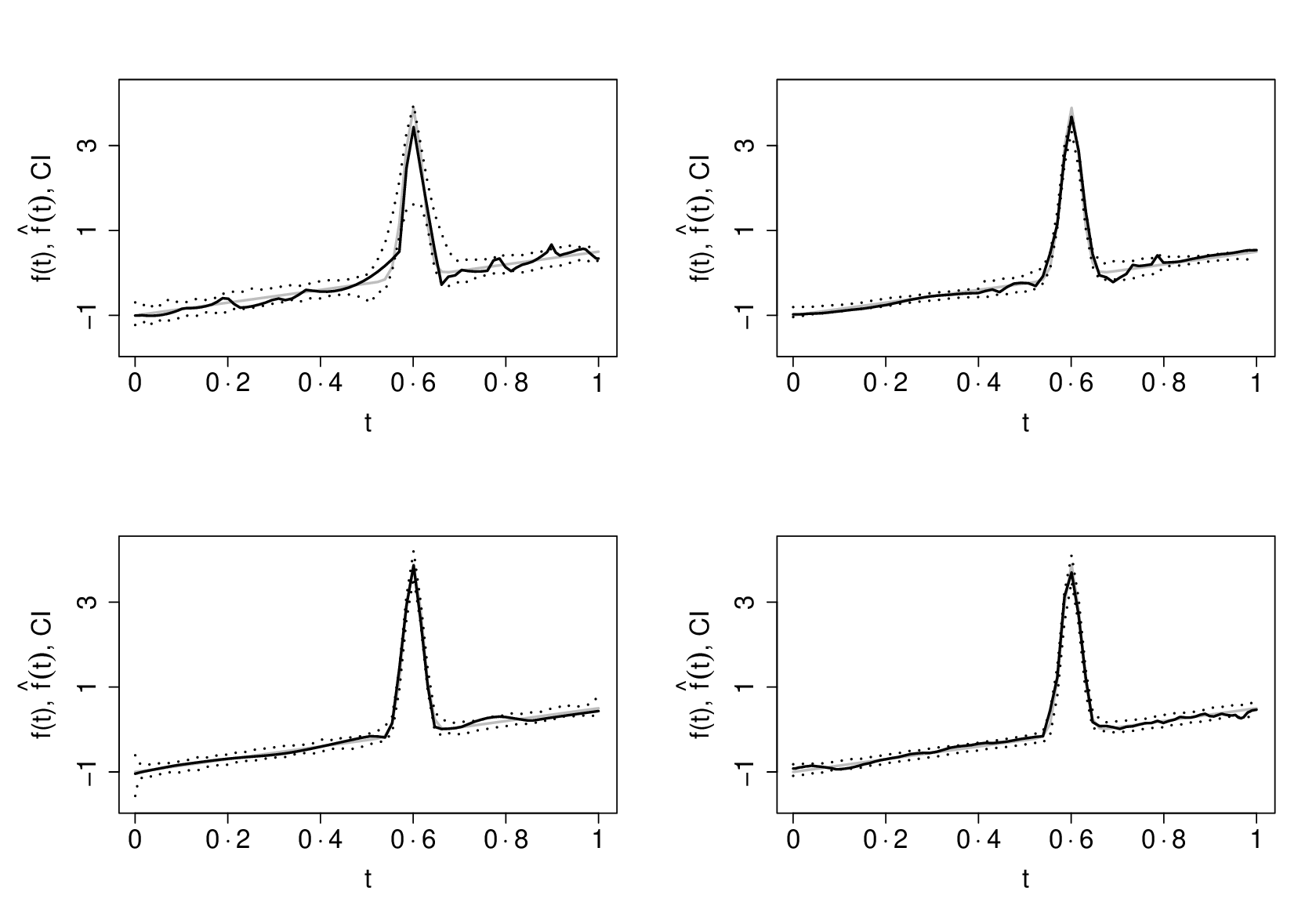,height=\linewidth,width=0.67\linewidth,clip=,angle=270}}
\caption{Estimates of the Mexican hat function for the data replicates with median integrated square errors.
The plotted curves are the true function (solid line), the spline estimate (solid line), and the
point-wise empirical 0.025 and 0.975 quantiles (dotted lines).
Top left: traditional smoothing spline estimate. Top right: estimate from the method in \cite{pintore_06}.
Bottom left: Loco-Spline estimate.
Bottom right: proposed adaptive smoothing spline estimate.
}
\label{fig:mexihat}
\end{figure}

Two well-known functions with
varying smoothness on the domain were considered under the model
$y_i=f(t_i)+\epsilon_i$ with $\epsilon_i\sim N(0,\sigma^2)$.
We used $n=200$ and $t_i=i/n$ $(i=1,\ldots,n)$ in all the simulations
and repeated each simulation on 100 randomly generated data replicates.
The integrated square error $\int_0^1 \{\hat{f}(t)-f_0(t)\}^2 dt$
and point-wise absolute errors at $t=0{\cdot}2,0{\cdot}4,0{\cdot}6,0{\cdot}8$ were
used to evaluate the performance of an estimate $\hat{f}$.
To visualize the comparison, we also selected for each example and each method
a data replicate with the median performance as follows. The function estimates from each method yielded 100 integrated square errors. After ranking them from the lowest to the highest, we chose the 50th integrated square error and its corresponding data replicate to represent the median performance.
We then plotted the function estimates from these selected data replicates in Fig.~\ref{fig:heav}-\ref{fig:mexihat} to compare the median estimation performances for different methods.
To assess variability in estimation, we also superimposed in these plots the point-wise empirical 0.025 and 0.975 quantiles of the 100 estimates.

We first consider data generated from the Heaviside function
$f(t)=5I_{[t>=0{\cdot}5]}$ with $\sigma=0{\cdot}7$. Based on the error summary statistics in Table~\ref{tab:sim}, all the adaptive methods outperform the traditional smoothing splines, with our method and that in \cite{pintore_06} displaying clear advantages in all the error measures.
Furthermore, our method had the smallest mean integrated square error. This advantage is better illustrated by the plots in Fig.~\ref{fig:heav}.
While the median estimates from all the three adaptive methods tracked the true function reasonably well, the Loco-Spline estimates show greater variability than the other two adaptive methods in estimating the flat parts of the Heaviside function. Further, our method does the best job in tracking down the jump. The estimate of \cite{pintore_06} can oscillate around the jump of the Heaviside function, probably because the equally-spaced
jump points for $\rho$ suggested in their paper sometimes have difficulty in characterizing the jump in the true function.
This echoes the finding in \cite{guo_10} that the jump locations of $\rho$ also need to be adaptive, a concept adopted in our method.

\begin{table}
\def~{\hphantom{0}}
\tbl{Comparison of integrated square errors and point-wise absolute errors for various estimates. Values, divided by 100, are empirical means and standard deviations (in brackets) based on
100 data replicates. }{%
\begin{tabular}{lccccc}
Method & ISE & PAE(0${\cdot}$2) & PAE(0${\cdot}$4) & PAE(0${\cdot}$6) & PAE(0${\cdot}$8)\\
\hline
& \multicolumn{5}{c}{Heaviside function}\\
SS & 18(7) & 15(11) & 17(14) & 16(14) & 16(12)\\
PSH & 5(2) & 6(5) & 6(5) & 7(5) & 7(5)\\
Loco & 7(3) & 10(8) & 13(12) & 11(10) & 12(12)\\
ADSS & 2(2)& 7(5) & 6(5)& 6(5) & 7(6)\\
\hline
& \multicolumn{5}{c}{Mexican hat function}\\
SS & 6${\cdot}$6(6${\cdot}$2) & 8(6) & 8(8) & 96(72) & 8(6)\\
PSH & 1${\cdot}$1(0${\cdot}$3) & 4(3) & 8(5) & 35(11) & 8(5)\\
Loco & 0${\cdot}$6(0${\cdot}$3) & 4(4) & 5(4) & 13(10) & 5(4)\\
ADSS & 0${\cdot}$6(0${\cdot}$2) & 4(3) & 4(3) & 15(10) & 6(4)\\
\hline
\end{tabular}}\label{tab:sim}
\begin{tabnote}
ISE, integrated square error; PAE, point-wise absolute error; SS, smoothing splines; PSH, splines in \cite{pintore_06}; Loco, Loco-Splines; ADSS,  adaptive smoothing splines in this paper
\end{tabnote}
\end{table}

The second example is the Mexican hat function
$f(t)=-1+1{\cdot}5t+0{\cdot}2\phi_{0{\cdot}02}(t-0{\cdot}6)$ with $\sigma=0{\cdot}25$, where
$\phi_{0{\cdot}02}(t-0{\cdot}6)$ is the density function of $N(0{\cdot}6,0{\cdot}02^2)$.
From Table~\ref{tab:sim} and Fig.~\ref{fig:mexihat}, the estimates from our method
and the Loco-Spline have competitive performance and both
outperform the traditional smoothing spline and those of \cite{pintore_06}. The estimates of \cite{pintore_06} again suffer close to the hat.

%

\begin{figure}
\centerline{\psfig{file=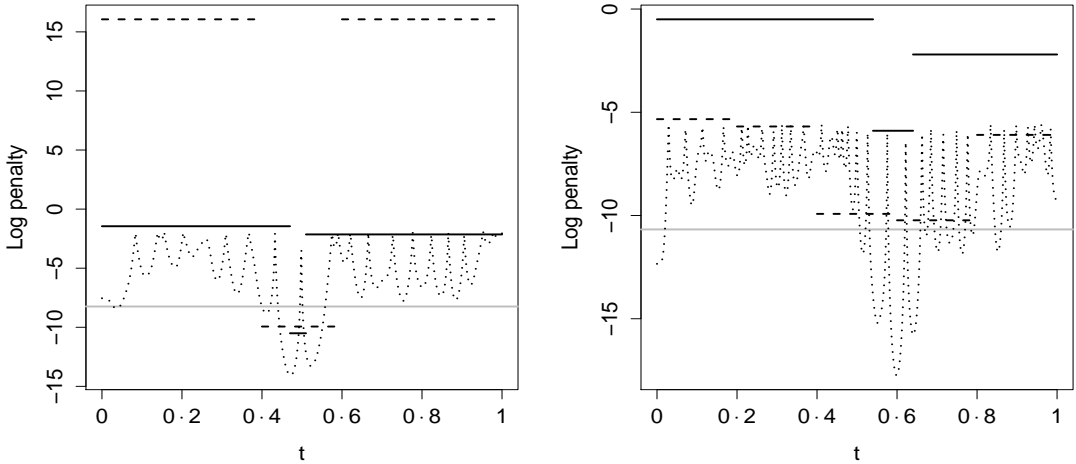,height=\linewidth,width=0.5\linewidth,clip=,angle=270}}
\caption{Estimated log penalties for simulation examples in Fig.~\ref{fig:heav}--\ref{fig:mexihat}.
The log penalties are for traditional smoothing splines (solid grey lines), the method in \cite{pintore_06} (dashed steps),
the Loco-Spline (dotted lines), and the proposed method (solid steps).
Left: Heaviside. Right: Mexican hat. 
}
\label{fig:penalty}
\end{figure}

For the estimates plotted in Fig.~\ref{fig:heav}--\ref{fig:mexihat}, we also
plot the estimated log penalties for all the methods 
in Figure~\ref{fig:penalty}.
In general, the penalty functions from the three adaptive methods track the smoothness changes in the underlying functions reasonably well.

\section{Application}

\begin{figure}
\centerline{\psfig{file=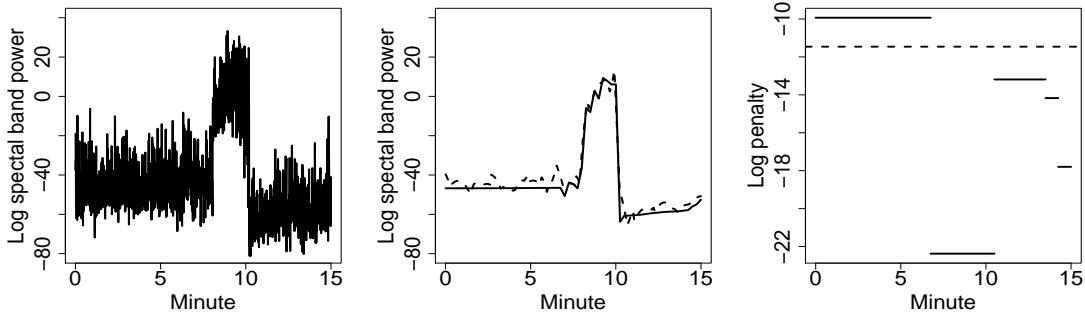,height=\linewidth,width=0.33\linewidth,clip=,angle=270}}
\caption{EEG data example.
Left: Raw log spectral band power. Center: Reconstructions from the traditional smoothing splines (dashed)
and the proposed adaptive smoothing splines (solid). Right: Estimated log penalties from the traditional smoothing splines and the proposed adaptive smoothing splines.
}
\label{fig:eeg}
\end{figure}

In this section, we apply the proposed adaptive smoothing splines
to an example on electroencephalograms of epilepsy patients
\citep{guo_10}. Previous research \citep{qin_09} has shown that the
low voltage frequency band 26-50Hz is important in characterizing
electroencephalograms and may help determine the spatial-temporal initiation of
seizure. The left panel of Figure~\ref{fig:eeg} shows the raw
time-varying log-spectral band power of 26-50Hz calculated every
half second for a 15-minute long intracranial electroencephalogram series. The
sampling rate was 200Hz and the seizure onset was at the 8th
minute \citep{litt_01}. The raw band powers are always very noisy
and need to be smoothed before further analysis. The middle panel
shows the reconstructions from traditional smoothing splines and the
proposed adaptive smoothing splines. We also tried the Loco-Spline
but the program exited due to a singular matrix error.

Traditional smoothing splines clearly under-smooth the pre- and post-seizure regions
and over-smooth the seizure period, because a single smoothing parameter is insufficient to capture the abrupt change before the onset of the seizure.
Our estimate smoothes out the noise on both ends but keeps the details before the onset of seizure. In particular, we see a fluctuation in power starting from a minute or so before the onset
of the seizure, which may be a meaningful predictor of seizure initiation. The band power then increases sharply
at the beginning of the seizure. Around the 10th minute at the end of the seizure,
the band power drops sharply to a level even lower than the
pre-seizure level, an indication of the suppression of neuronal activities after seizure. Afterwards, the band power starts
to regain. But it still fails to reach the pre-seizure level even at the end of the 15th minute.
These findings concur with those in \cite{guo_10}.

The proposed method took less than 10 minutes for the whole analysis, compared with 40-50 minutes for the method in \cite{guo_10}. This is not surprising, since the latter not only needs a dense grid search to locate the jump points but also lacks good initial step sizes.

\section*{Acknowledgements}

We are grateful to two referees and Associate Editor for
constructive and insightful comments. We are also thankful to
Wensheng Guo and Ziyue Liu for providing the electroencephalogram data, and Howard
Bondell for help with the Loco-spline program. Xiao Wang's research is supported by US NSF grants CMMI-1030246 and DMS-1042967 and Jinglai Shen's research is supported by US NSF grants CMMI-1030804 and DMS-1042916.

\section*{Supplementary material}

Supplementary Material available at {\it Biometrika} online includes the proofs of Theorems 1-3 and Corollary 1, and the detailed derivation of the Green's function.

\appendix
%
%
\appendixone
\section*{Appendix}


In this appendix, we provide  outline proofs of Theorems~\ref{thm:char} and \ref{thm:asy-rep}.
For the full proofs of these two theorems and Corollary~\ref{coro:asy}, we refer the readers to the Supplementary Material.

{\it Outline Proof of Theorem \ref{thm:char}}.
For any $f, g\in W^m_2$ and $\delta \in \mathbb R$,
\begin{equation}\label{equ:temp1}
\psi(f+\delta g) - \psi(f) = 2\delta \psi_1(f, g) + \delta^2\Big[\int_0^1 g^2(t)d\omega_n(t) + \lambda \int_0^1 \rho(t)\{g^{(m)}(t)\}^2dt\Big],
\end{equation}
where
\begin{equation}\label{equ:psi1}
\psi_1(f, g) = \int_0^1 \sigma^{-2}(t)\{f(t) - h(t)\}g(t) d\omega_n(t) +
\lambda \int_0^1 \rho(t) f^{(m)}(t)g^{(m)}(t)dt.
\end{equation}

\begin{lemma}\label{lem:psi1} The function $f\in W^m_2$ minimizes $\psi(f)$ in (\ref{equ:obj}) if and only if $\psi_1(f, g)=0$ for all $g\in W_2^m$.
\end{lemma}

Let $g(t) = t^k (k=0, \ldots, m-1)$ in \eqref{equ:psi1}. An application of Lemma A\ref{lem:psi1} shows that if $f$ minimizes
$\psi(f)$, then
$$\int_0^1 \sigma^{-2}(t)\{f(t) - h(t)\}~t^k d\omega_n(t) = 0 ~~~ (k=0, 1, \ldots, m-1).$$
We first have $$\check l_1(f, 1) - \check l_1(h, 1) = \int_0^1 \sigma^{-2}(t)\{f(t) - h(t)\} d\omega_n(t) = 0.$$
Further,
$$\check l_2(f, 1) - \check l_2(h, 1)  = \int_0^1 \int_0^s\sigma^{-2}(t)\{f(t) - h(t)\} d\omega_n(t)ds= \int_0^1  \sigma^{-2}(t)\{f(t) - h(t)\} ~t ~d\omega_n(t)
= 0.$$
Similarly, $\check l_k(f, 1) = \check l_k(h, 1)$ for $k=1, \ldots, m$.

\begin{lemma} If $f\in W_2^m$ satisfies $\check l_k(f, 1) = \check l_k(h, 1)$, $k=1, \ldots, m$, then for all $g\in W_2^m$,
\begin{equation}\label{equ:psi2-1}
\psi_1(f, g) = \int_0^1 \psi_2(f) ~g^{(m)}(t) dt,
\end{equation}
where
\begin{equation}\label{equ:psi2-2}
\psi_2(f) = \lambda ~\rho(t)~ f^{(m)}(t) + (-1)^m~ \{\check l_m(f, t) - \check l_m(h, t)\}.
\end{equation}
\end{lemma}

Let $B^{+} = \{t\in [0, 1] :  \psi_2(f) > 0\}$ and $B^{-} = \{t \in [0, 1] :  \psi_2(f) < 0\}$. Define $g_+^{(m)}(t) = -I_{B^+}(t)$ and $g_{-}^{(m)}(t) = I_{B^{-}}(t)$, where $I$ is the indicator function.
Since $\psi_1(f, g) =0$ for all $g\in W_2^m$,  we have $\psi_1(f, g_+)<0$ and $\psi_1(f, g_-)<0$, unless $B^+$ and $B^-$ are of measure zero.
This shows that $\psi_2(f)=0$ almost everywhere.


{\it Outline Proof of Theorem \ref{thm:asy-rep}}. It follows from \eqref{equ:rep0} that
$r^{-1}(t) \hat f(t)= V_1(t) + V_2(t) + V_3(t)+V_4(t)$, where
\begin{align*}
V_1(t) &= {d^m\over dt^m}\int_0^1 P(t, s) l_m(f_0, s)ds, ~~V_2(t) = {d^m\over dt^m}\int_0^1 P(t, s)\{\check l_m(h, s) - \check l_m(f_0, s)\}ds,\\
V_3(t) &=  {d^m\over dt^m} \int_0^1 P(t, s)\{l_m(\hat f - f_0, s)-\check l_m(\hat f - f_0, s)\}ds,~~ V_4(t) = \sum_{k=1}^{2m} a_k C_k^{(m)}(t).
\end{align*}
Let $\bar f$ minimize the functional
$$\int_0^1 r^{-1}(s)\{f(s) - f_0(s)\}^2 ds + \lambda \int_0^1 \rho(t) f^{(m)}(s)^2ds.$$
Similar to Theorem \ref{thm:char}, we have
\begin{equation}\label{equ:asy-mean}
(-1)^m \lambda \rho(t) \bar f^{(m)}(t) + l_m(\bar f, t) = l_m(f_0, t),\end{equation} and
\begin{equation}\label{equ:asy-mean2} l_m(\bar f, t) = \int_0^1 P(t, s) l_m(f_0, s)ds.\end{equation}
Hence, $V_1(t) = r^{-1}(t) \bar f(t)$. Taking the $m$th
derivative of both sides of \eqref{equ:asy-mean},  we get
$$(-1)^m \lambda \{\rho(t)\bar f^{(m)}(t)\}^{(m)} + r^{-1}(t) \bar f(t) =r^{-1}(t)f_0(t).$$
Recall that $f_0$ is $2m$ times continuously differentiable and $\beta = \lambda^{-1/(2m)}$. Combining this with \eqref{equ:asy-mean2}, it is easy to show that  $\bar
f^{(k)}(t) \rightarrow f_0^{(k)}(t)$ as $\beta\rightarrow
\infty$ for $k=1, \ldots, 2m$. Therefore,
$$V_1(t) = r^{-1}(t) f_0(t) + (-1)^{m-1} \lambda \{\rho(t) f_0^{(m)}(t)\}^{(m)} +
o(\lambda).$$

\begin{proposition}\label{prop:change} Assume that a function $\tilde J(t, s)$ satisfies
$(-1)^m{\partial^m\over \partial s^m} \tilde J(t, s) = {\partial^m\over \partial t^m}P(t, s), ~~ t, s \in [0, 1].$
Then
$\tilde J(t, s)  + \sum_{k=0}^{m-1} (-1)^{k} \zeta_{k+1}(s)\tilde J_{k}(t) =  
( r(s)/r(t) ) J(t, s),$
where
$$\zeta_k(s) = \int_{s}^1\cdots\int_{s_{k-3}}^1\int_{s_{k-2}}^1 d s_{k-1} ds_{k-2}\cdots ds_1,~~~~\tilde J_k(t) = {\partial^{k}\over \partial s^{k}}\tilde J(t, s)\mid_{s=1},$$
and $J(t, s)$ is the Green's function for
\begin{equation}\label{equ:ode0}
(-1)^m \lambda r(t) \{\rho(t)
u^{(m)}(t)\}^{(m)} + u(t) =0. \end{equation}
\end{proposition}



By applying Proposition A\ref{prop:change}, we have, for any $t\in
(0, 1)$,
\begin{eqnarray*}
V_2(t) &=& \int_0^1 (-1)^m {\partial^m\over \partial s^m}\tilde J(t, s) \check l_m(h-f_0, s)ds\\
&=& \int_0^1 \tilde J(t, s) d\{\check l_1(h - f_0, s)\} + (-1)^m \sum_{k=1}^{m-1} (-1)^{k-1} \tilde J_{m-k}(t) \check l_{m-k+1}(h - f_0, 1)
\\
&=& {1\over n}\sum_{i=1}^n {r(t_i)\over r(t)} J(t, t_i)\sigma^{-1}(t_i) \epsilon_i + \mbox{ higher order terms. }
\end{eqnarray*}
\cite{eggermont_06} established the uniform error bounds for
regular smoothing splines.  We adopt the same approach as in
\cite{eggermont_06} for adaptive smoothing splines; the details
are omitted here. For $\lambda \ll
(n^{-1}\log n)^{2m/(1+4m)}$, we obtain
$$\|\hat f - f_0\| =  O\left[\Big\{{\max \big(\log{1\over \lambda}, \log\log n \big) \over n\lambda^{1/(2m)}}~\Big\}^{1/2}\right].$$
Therefore, $\|V_3\| \le O(\beta^m) D_n \|\hat f - f_0\|$. Finally,
it is shown in detail in the Supplementary Material
that $\|V_4\|$ is of order $O(\beta^m) \exp[-\beta
Q_\beta(t)\{Q_\beta(1)- Q_\beta(t)\}]$, and thus a
negligible term in the asymptotic expansion of $r^{-1}(t)\hat{f}(t)$. This completes the representation for $\hat f$.

\bibliographystyle{biometrika}

\begin{thebibliography}{99}





\bibitem[\protect\astroncite{Abramovich \&\ Grinshtein}{1999}]{abramovich_99}
{\sc Abramovich, F. \& Grinshtein, V.} (1999). Derivation of an equivalent
kernel for general spline smoothing: a systematic approach. {\it
Bernoulli} {\bf 5}, 359--79.



\bibitem[\protect\astroncite{Berlinet \&\ Thomas-Agnan}{2004}]{berlinet_04}
{\sc Berlinet, A. \& Thomas-Agnan, C.} (2004). {\em Reproducing
Kernel Hilbert Spaces in Probability and Statistics}, Kluwer: Wiley.

\bibitem[\protect\astroncite{Cai \&\ Wang}{2008}]{cai_08}
{\sc Cai, T. \& Wang, L.} (2008). Adaptive variance function estimation in heteroscedastic nonparametric regression. {\it Ann. Statist.} {\bf 36}, 2025--54.

%
\bibitem[\protect\astroncite{Coddington \&\ Levinson}{1955}]{Coddington_book55}
{\sc Coddington, E. A.  \& Levinson, N.} (1955). {\em Theory of Ordinary
Differential Equations}. New York: McGraw-Hill.
%

\bibitem[\protect\astroncite{Cox}{1983}]{cox_83}
{\sc Cox, D. D.} (1983). Asymptotics of M-type smoothing splines. {\it Ann. Statist.} {\bf 11}, 530--51.

\bibitem[\protect\astroncite{Craven \& Wahba}{1979}]{craven_79}
{\sc Craven, P. \& Wahba, G.} (1979). Smoothing noisy data with spline functions: estimating the correct degree
of smoothing by the method of generalized cross-validation. {\it Numer. Math.} {\bf 31}, 377--403.

\bibitem[\protect\astroncite{DiMatteo et al.}{2001}]{DiMatteo_01}
{\sc DiMatteo, I., Genovese, C. R. \& Kass, R. E.} (2001). Bayesian
curve-fitting with free-knot splines. {\it Biometrika} {\bf 88},
1055--71.

\bibitem[\protect\astroncite{Donoho \&\ Johnstone}{1994}]{donoho_94}
{\sc Donoho, D. L. \& Johnstone, I. M.} (1994). Ideal spatial adaptation by wavelet shrinkage. {\it Biometrika} {\bf 81}, 425--55.

\bibitem[\protect\astroncite{Donoho \&\ Johnstone}{1995}]{donoho_95}
{\sc Donoho, D. L. \& Johnstone, I. M.} (1995). Adaptive to unknown smoothness via wavelet shrinkage. {\it J. Ameri. Statist. Assoc.} {\bf 90}, 1200--24.

\bibitem[\protect\astroncite{Donoho \&\ Johnstone}{1998}]{donoho_98}
{\sc Donoho, D. L. \& Johnstone, I. M.}(1998). Minimax estimation via wavelet shrinkage. {\it Ann. Statist.} {\bf 26}, 879--921.

\bibitem[\protect\astroncite{Eggermont \&\ LaRiccia}{2006}]{eggermont_06}
{\sc Eggermont, P. P. B.  \&  LaRiccia, V. N.}(2006).  Uniform error bounds
for smoothing splines. {\it IMS Lecture Notes-Monograph Series:
High Dimensional Probability} {\bf 51}, 220--37.

\bibitem[\protect\astroncite{Eggermont \&\ LaRiccia}{2009}]{eggermont_09}
{\sc Eggermont, P. P. B.  \&  LaRiccia, V. N.} (2009). {\it Maximum Penalized
Likelihood Estimation}. Volume II: Regression. New York: Springer.

\bibitem[\protect\astroncite{Eubank}{1999}]{eubank_99} {\sc Eubank, R. L.} (1999). {\it Nonparametric Regression and Spline Smoothing}. New York: Marcek Dekker.

\bibitem[\protect\astroncite{Fan \&\ Gijbels}{1996}]{fan_96}
{\sc Fan, J. \& Gijbels, I.} (1996). {\em Local Polynomial Modelling and Its Applications}.
Chapman \& Hall, Boca Raton.

\bibitem[\protect\astroncite{Fan \&\ Yao}{1998}]{fan_98}
{\sc Fan, J. \& Yao, Q.} (1998). Efficient estimation of conditional
variance functions in stochastic regression. {\it Biometrika}
{\bf 85}, 645--60.

\bibitem[\protect\astroncite{Friedman \&\ Silverman}{1989}]{friedman_89}
{\sc Friedman, J. \& Silverman, B. W.} (1989). Flexible parsimonious smoothing and additive modeling (with discussion). {\it Technometrics} {\bf 31}, 3--39.

\bibitem[\protect\astroncite{Gu}{2002}]{gu_02} {\sc Gu, C.} (2002). {\it Smoothing Spline ANOVA Models}. New York: Springer.

\bibitem[\protect\astroncite{Hall \&\ Carroll}{1989}]{hall_89}
{\sc Hall, P. \& Carroll, R. J.}(1989). Variance function estimation in
regression: the effect of estimating the mean. {\it J. Roy.
Statist. Soc. Ser. B} {\bf 51}, 3--14.

\bibitem[\protect\astroncite{Hansen \&\ Kooperberg}{2002}]{hansen_02}
{\sc Hansen, M. H. \& Kooperberg, C.}(2002). Spline adaption in extened
linear models. {\it Statist. Sci.} {\bf 17}, 2--51.

\bibitem[\protect\astroncite{Kimeldorf \&\ Wahba}{1971}]{kim_71}
{\sc Kimeldorf, G. S. \& Wahba, G.}(1971). Some results on Tchebycheffian spline functions. {\it J Math.
Anal. Applic.} {\bf 33}, 82--95.

\bibitem[\protect\astroncite{Kohn \&\ Ansley}{1983}]{kohn_83}
{\sc Kohn, R. \& Ansley, C. F.} (1983). On the smoothness properties of the best linear unbiased estimate
of a stochastic process observed with noise. {\it Ann. Statist.} {\bf  11}, 1011--17.

\bibitem[\protect\astroncite{Litt et al}{2001}]{litt_01}
{\sc Litt, B., Esteller R., Echauz, J., D'Alessandro, M.,  Shor, R.,
Henry, T., Pennell, P., Epstein, C., Bakay, R., Dichter, M. \&
Vachtsevanose, G.}(2001). Epileptic seizures may begin hours in advance of
clinical onset: a report of five patients. {\it Neuron} {\bf 30},
51--64.

\bibitem[\protect\astroncite{Liu \&\ Guo}{2010}]{guo_10}
{\sc Liu, Z. \&  Guo, W.} (2010). Data driven adaptive spline smoothing.
{\it Statistica Sinica} {\bf 20}, 1143--63.

\bibitem[\protect\astroncite{Luo \&\ Wahba}{1997}]{luo_97}
{\sc Luo, Z. \& Wahba, G.} (1997). Hybrid adaptive splines. {\it J. Am. Statist. Assoc.} {\bf 92}, 107--16.

\bibitem[\protect\astroncite{Mao \&\ Zhao}{2003}]{mao_03}
{\sc Mao, W. \& Zhao, L.} (2003). Free-knot polynomial splines with confidence
intervals. {\it J. R. Statist. Soc. B} {\bf 65}, 901--19.


\bibitem[\protect\astroncite{Messer}{1991}]{messer_91} {\sc Messer, K.} (1991). A comparison of a spline estimate to
its equivelent kernel estimate. {\it Ann. Statist.} {\bf
19}, 817--29.


\bibitem[\protect\astroncite{M\"uller \&\ Stadtm\"uller}{1987}]{muller_87}
{\sc M\"uller, H.-G.  \& Stadtm\"uller, U.} (1987). Variable bandwidth kernel estimators
of regression curves. {\it Ann. of Statist.} {\bf 15}, 282--301.

\bibitem[\protect\astroncite{Nussbaum}{1985}]{nussbaum_85}
{\sc Nussbaum, M.} (1985). Spline smoothing in regression models
and asymptotic efficiency in $L_2$. {\it Ann. Statist.}
{\bf 13}, 984--97.

\bibitem[\protect\astroncite{Nychka}{1995}]{nychka_95} {\sc Nychka, D.}(1995). Splines as local smoothers. {\it Ann.
Statist.} {\bf 23}, 1175--97.

\bibitem[\protect\astroncite{Pinsker}{1980}]{pinsker_80}
{\sc Pinsker, M. S.} (1980). Optimal filtering of square
integrable signals in Gaussian white noise. {\it Problems Inform.
Transmission} {\bf 16}, 120--33.



\bibitem[\protect\astroncite{Pintore et al.}{2006}]{pintore_06}
{\sc Pintore, A., Speckman, P. \& Holmes,  C. C.} (2006). Spatially adaptive
smoothing splines. {\it Biometrika} {\bf 93}, 113--25.

\bibitem[\protect\astroncite{Qin et al.}{2009}]{qin_09}
{\sc Qin, L. , Guo, W. \& Litt, B.} (2009). A time-frequency functional model for locally stationary
time series data. {\it J. Comput. Graph. Statist.}, {\bf  18}, 675--93.

\bibitem[\protect\astroncite{Rice \&\ Rosenblatt}{1983}]{rice_83} {\sc Rice, J. \& Rosenblatt, M.}(1983). Smoothing
splines: regression, derivatives and deconvolution. {\it Ann. Statist.} {\bf 11},  141--56.

\bibitem[\protect\astroncite{Ruppert \&\ Carroll}{2000}]{ruppert_00}
{\sc Ruppert, D. \& Carroll, R. J.} (2000). Spatially-adaptive penalties for spline
fitting. {\it Aust. N. Z. J. Statist.} {\bf 42}, 205--23.

\bibitem[\protect\astroncite{Ruppert et al.}{1997}]{ruppert_97}
{\sc Ruppert, D., Wand, M. P.,  Holst, U.  \& H$\mathrm{\ddot{o}}$sjer, O.} (1997). Local polynomial variance
function estimation. {\it Technometrics} {\bf 39}, 262--73.

\bibitem[\protect\astroncite{Silverman}{1984}]{silverman_84} {\sc Silverman, B. W.}(1984). Spline smoothing: the equivalent variable kernel
method. {\it Ann.  Statist.} {\bf 12}, 898--916.

\bibitem[\protect\astroncite{Silverman}{1985}]{silverman_85} {\sc Silverman, B. W.} (1985). Some aspects of the spline smoothing approch to nonparametric
curve fitting (with discussion). {\it J.  R. Statist. Soc. B} {\bf 47}, 1--52.

\bibitem[\protect\astroncite{Smith \&\ Kohn}{1996}]{smith_96}
{\sc Smith, M. \& Kohn, R.} (1996). Nonparametric regression using Bayesian variable selection. {\it J. Economet.}
{\bf 75},  317--44.

\bibitem[\protect\astroncite{Speckman}{1981}]{speckman_81}
{\sc Speckman, P. L.} (1981). The asymptotic integrated mean squared error for smoothing noisy data by splines. {\it Technical Report}, University of Oregon.

\bibitem[\protect\astroncite{Stone}{1982}]{stone_82}
{\sc Stone, C. J.} (1982). Optimal rate of convergence for nonparametric regression. {\it
Ann. Statist.} {\bf 10},  1040--53.


\bibitem[\protect\astroncite{Stone et al.}{1997}]{stone_97}
{\sc Stone, C. J.,  Hansen, M.,  Kooperberg, C. \& Truong, Y. K.} (1997). Polynomial splines and their tensor products in extended linear models. {\it Ann. Statis.} {\bf 25},  1371--425.


\bibitem[\protect\astroncite{Storlie et al.}{2010}]{Storlie_09JCGS} {\sc Storlie, C. B.,  Bondell, H. D.  \& Reich, B. J.}(2010). A
locally adaptive penalty for estimation of functions with varying
roughness. {\it J. Comp. Grap.
Statist.} {\bf 19}, 569--89.

\bibitem[\protect\astroncite{Wahba}{1985}]{wahba_85}
{\sc Wahba, G.}(1985). A comparison of GCV and GML for choosing the smoothing parameter in the generalized
spline smoothing problem. {\it Ann. Statist.}  {\bf 13}, 1378--402.

\bibitem[\protect\astroncite{Wahba}{1990}]{wahba_90}
{\sc Wahba, G.}(1990). {\it Spline Models for Observation Data}. Society
for Industrial and Applied Mathematics.

\bibitem[\protect\astroncite{Wahba}{1995}]{wahba_95}
{\sc Wahba, G.}(1995). Discussion of `Wavelet shrinkage: asymptopia?' by D. L. Donoho, I. M. Johnstone, G. Kerkyacharian \& D. Picard. {\it J. R. Statist. Soc. B} {\bf 57}, 360--1.

\bibitem[\protect\astroncite{Wang et al.}{2010}]{wang_10}
{\sc Wang, X., Shen, J. \& Ruppert, D.}(2010). Local asymptotics of
$P$-splines. {\it Electr. J. Statist.} {\bf  5},
1--17.

\bibitem[\protect\astroncite{Wood et al.}{2002}]{wood_02}
{\sc Wood, S. A., Jiang, W. \& Tanner, M.}(2002). Bayesian mixture of splines for spatially adaptive nonparametric
regression. {\it Biometrika} {\bf  89}, 513--21.

\end{thebibliography}

\end{document}